\documentclass[preprint,12pt]{elsarticle}

\usepackage{hyperref,amsmath,amsfonts,mathrsfs}

\usepackage{graphicx}
\usepackage{amsfonts}
\usepackage{amsbsy}
\usepackage{amssymb}
\usepackage{amsmath,tikz,pgfplots}
\usetikzlibrary{decorations.pathmorphing}
\usepackage{tikz}
\numberwithin{equation}{section}
\usetikzlibrary{arrows}
\usepackage{graphics}
\usepackage{graphicx}
\usepackage{amsfonts}
\usepackage{amsbsy}
\usepackage{amssymb}
\usepackage{amsmath,tikz,pgfplots}
\usetikzlibrary{decorations.pathmorphing}
\usepackage{graphics}
\def\bpsp{\begin{pspicture}}
\def\epsp{\end{pspicture}}

\newtheorem{theorem}{Theorem}[section]
\newtheorem{remark}[theorem]{Remark}
\newtheorem{example}[theorem]{Example}
\newtheorem{lemma}[theorem]{Lemma}
\newtheorem{corollary}[theorem]{Corollary}
\newtheorem{definition}[theorem]{Definition}
\newtheorem{proposition}[theorem]{Proposition}
\newtheorem{note}{Note}
\newtheorem{case}{Case}
\newtheorem{conjecture}{Conjecture}
\newtheorem{question}{Question}

\newcommand{\bea}{\begin{eqnarray}}
\newcommand{\eea}{\end{eqnarray}}
\newcommand{\beq}{\begin{eqnarray*}}
\newcommand{\eeq}{\end{eqnarray*}}

\def\m4{\mbox{\rm ~(mod $4$)}}

\def \bd{\begin{definition}}
\def \ed{\end{definition}}
\def \bqu{\begin{question}}
\def \equ{\end{question}}
\def \bcc{\begin{conjecture}}
\def \ecc{\end{conjecture}}
\def \bt{\begin{theorem}}
\def \et{\end{theorem}}
\def \bl{\begin{lemma}}
\def \el{\end{lemma}}
\def \bc{\begin{corollary}}
\def \ec{\end{corollary}}
\def \be{\begin{equation}}
\def \ee{\end{equation}}
\def \ben{\begin{enumerate}}
\def \een{\end{enumerate}}
\def \ba{\begin{array}}
\def \ea{\end{array}}
\def \bp{\begin{proposition}}
\def \ep{\end{proposition}}
\def \bx{\begin{example}}
\def \ex{\end{example}}
\def \br{\begin{remark}}
\def \er{\end{remark}}
\def \bdsc{\begin{description}}
\def \edsc{\end{description}}

\def \bn{\begin{case}}
\def \en{\end{case}}
\def \bnt{\begin{note}}
\def \ent{\end{note}}
\def\1{1\!\!1}

\def\mm2{\mbox{\rm ~(mod $2$)}}
\def\m4{\mbox{\rm ~(mod $4$)}}

\def\qed{\nolinebreak\hfill\rule{.2cm}{.2cm}\par\addvspace{.5cm}}

\def\m{\mu}

\def\1{\textbf{1}}
\def\0{\textbf{0}}

\journal{ABC}

\begin{document}

\begin{frontmatter}



\title{Variation in $\alpha$ trace norm of a digraph  by deletion of a vertex or an arc and its applications}
\author{Mushtaq A. Bhat$^{a}$}
\author{Peer Abdul Manan$^{b}$}

\address{Department of  Mathematics,  National Institute of Technology, Srinagar-190006, India}
\address {$^{a}$mushtaqab@nitsri.ac.in;~~~$^{b}$\text{mananab214@gmail.com}}

\begin{abstract}
 Let $D$ be a digraph of order $n$ with adjacency matrix $A(D)$. For $\alpha\in[0,1)$, the $A_{\alpha}$ matrix of $D$ is defined as $A_{\alpha}(D)=\alpha {\Delta}^{+}(D)+(1-\alpha)A(D)$, where ${\Delta}^{+}(D)=\mbox{diag}~(d_1^{+},d_2^{+},\dots,d_n^{+})$ is the diagonal matrix of vertex outdegrees of $D$. Let $\sigma_{1\alpha}(D),\sigma_{2\alpha}(D),\dots,\sigma_{n\alpha}(D)$ be the singular values of $A_{\alpha}(D)$. Then the trace norm of $A_{\alpha}(D)$, which we call $\alpha$ trace norm of $D$, is defined as $\|A_{\alpha}(D)\|_*=\sum_{i=1}^{n}\sigma_{i\alpha}(D)$. In this paper, we study the variation in $\alpha$ trace norm of a digraph when a vertex or an arc is deleted. As an application of these results, we characterize oriented trees and unicyclic digraphs with maximum $\alpha$ trace norm.
\end{abstract}

\vskip 0.2 true cm

\begin{keyword} $\alpha$ Singular values, $\alpha$ Trace norm, Oriented tree, Unicyclic digraph.

\vskip 0.2 true cm


$MSC$: 05C20, 05C50 

\end{keyword}

\end{frontmatter}

\section{\bf Introduction}
Let $D(\mathcal{V},\mathcal{A})$ be a digraph with vertex set $\mathcal{V}=\{v_1,v_2,\dots,v_n\}$ and arc set $\mathcal{A}$. Our digraphs are simple i.e., loop free and free of parallel arcs. In a digraph $D$, If there is an arc from vertex $u$ to vertex $v$, we denote this by $uv$. An oriented graph is a digraph with no pair of symmetric arcs.  
 A graph $G$ can be regarded as a symmetric digraph $\overleftrightarrow {G}$ obtained by replacing each edge $e$ of $G$ by a pair of symmetric arcs. The set of vertices $\{w\in \mathcal{V}: uw\in \mathcal{A}\}$ is called the outneighbour of $u$ and we denote it by $N^{+}(u)$. In a digraph, the set of vertices $\{w\in \mathcal{V}: wu\in \mathcal{A}\}$ is called the inneighbour of $u$ and we denote this by $N^{-}(u)$. The cardinality of the set $N^{+}(u)$ is called outdegree of $u$ and we denote it by $d_u^+$.  The cardinality of the set $N^{-}(u)$ is called indegree of $u$ and we denote it by $d_u^-$.\\
\indent Let $D$ be a digraph with vertex set $\mathcal{V}=\{v_1,v_2,\dots,v_n\}$. Then the adjacency matrix $A(D)=(a_{ij})$ of $D$ is a square matrix of order $n$ with $a_{ij}=1$ if there is an arc from vertex  $v_i$ to vertex $v_j$ and zero, otherwise. Let ${\Delta}^+={\Delta}^+(D)=\mbox{diag}~(d^+_1,d^+_2,d^+_3,\dots,d^+_n)$, where $d^+_i=d^+_{v_i}$, be the diagonal matrix of vertex outdegrees of $D$. Then the Laplacian and signless Laplacian matrices of $D$ are respectively defined as $L(D)={\Delta}^{+}-A(D)$ and $Q(D)={\Delta}^{+}+A(D)$. For $\alpha\in[0, 1]$, Nikiforov \cite{n2} defined the alpha matrix $A_{\alpha}(G)$ of a graph $G$ as a common extension of adjacency matrix and signless Laplacian matrix as $A_{\alpha}(G)=\alpha D(G)+(1-\alpha)A(G)$, where $D(G)$ is diagonal matrix of vertex degrees of graph $G$. Liu et al. \cite{lwc}, defined the $A_{\alpha}$ matrix for digraphs and studied $\alpha$-spectral radius of digraphs. For $\alpha\in [0, 1)$, the alpha matrix of a digraph is defined as $A_{\alpha}(D)=\alpha{\Delta}^+(D)+(1-\alpha)A(D)$. For $\alpha=0$, we see $A_0(D)=A(D)$ and for $\alpha=\frac{1}{2}$, $A_{\frac{1}{2}}(D)=\frac{1}{2}Q(D)$. It is clear that $A_{\alpha}(D)$ is a common extension of adjacency matrix $A=A(D)$ and signless Laplacian matrix $Q(D)$ of a digraph $D$. Very recently, Yang et al. \cite{ybw} studied the spectral moments of $A_{\alpha}$ matrix of digraphs.\\
\indent Let $G$ be an undirected graph of order $n$ and with eigenvalues $\lambda_1\ge\lambda_2\ge\lambda_3\ge\dots\ge\lambda_n$. Then the energy of $G$ is defined as $E(G)=\sum_{i=1}^{n}|\lambda_i|$. This concept was given by Gutman (1978). For details related to graph energy see \cite{lsg}. The concept of energy was extended to digraphs by Pe\~{n}a and Rada \cite{pr} and they defined the energy of a digraph $D$ as $E(D)=\sum_{i=1}^{n}|\Re z_i|$, where $z_1,z_2,\dots,z_n$ are eigenvalues of $D$, possibly complex and $\Re z_i$ denote the real part of complex number $z_i$.\\
Given a complex matrix $M\in M_n(\mathbb{C})$, its trace norm is defined as 
\begin{equation*}
\|{M} \|_* =\sum\limits_{i=1}^{n}\sigma_i(M),
\end{equation*}
where $\sigma_1(M)\ge\sigma_2(M)\ge \sigma_3(M)\ge \dots\ge \sigma_n(M)$ are the singular values of $M$ i.e., the positive square roots of eigenvalues of $MM^{*}$.\\
Throughout the paper, we call singular values of $A_{\alpha}(D)$ as the $\alpha$ singular values of $D$ and trace norm of $A_{\alpha}(D)$ as the $\alpha$ trace norm of $D$ and we will denote the trace norm of $A_{\alpha}(D)$ by $\|D_{\alpha}\|_{*}$. If a digraph $D$ has $k$ distinct singular values $\sigma_{1\alpha}, \sigma_{2\alpha}, \dots, \sigma_{k\alpha}$, with their respective multiplicities $m_1, m_2,\dots,m_k$, then we write the set of singular values as $$\{{\sigma^{[m_1]}_{1\alpha}}, {\sigma^{[m_2]}_{2\alpha}},\dots,{\sigma^{[m_k]}_{k\alpha}}\}.$$
If $M=A(G)$, the adjacency matrix of a graph $G$, then $\sigma_i(M)=|\lambda_i(G)|$ and so trace norm coincides with the energy of a graph. Trace norm of a matrix is also known as Nikiforov energy of a matrix \cite{ar,n1}. So, graph energy extends to digraphs via trace norm as well. By $\mathcal{T}(n)$ and $\mathcal{U}(n)$, we respectively denote the set of oriented trees on $n$ vertices and set of unicyclic digraphs on $n$ vertices. Agudelo and Rada \cite{ar} obtained lower bounds for trace norm of digraphs and characterized the extremal digraphs. Agudelo et al. \cite{apr}, studied variation in trace norm of an oriented graph by vertex (leaf) deletion and determined oriented trees in $\mathcal{T}(n)$ attaining maximum and minimum trace norm. Monsalve and Rada \cite{mr1}, obtained oriented bipartite graphs with minimum trace norm. Monsalve et al. \cite{mr2} studied variation in trace norm by arc deletion and obtained oriented graphs with maximum and minimum trace norm in the class of oriented bicyclic graphs. Recently, Rather et al. \cite{rgw} studied trace norm for eccentricity matrices of graphs. Bhat and Manan \cite{bm} studied trace norm for $A_{\alpha}$ matrices of digraphs.\\
Let $\overrightarrow{P_n}$ denote the directed path with $n$ vertices $v_1,v_2,\dots,v_n$ and arc set $\mathcal{A}=\{v_1v_2,v_2v_3,\dots,\\v_{n-1}v_n\}$ and let $\overrightarrow{C_n}$ denote a directed cycle with $n$ vertices  $v_1,v_2,\dots,v_n$ and arc set $\mathcal{A}=\{v_1v_2,v_2v_3,\dots,v_{n-1}v_n, v_nv_1\}$. We first find the variation in $\alpha $ trace norm of a digraph by deleting an arc. Using this results, we determine variation in $\alpha$ trace norm of a digraph by deleting a vertex. Using these results, we determine the unique oriented tree in $\mathcal{T}(n)$ and unique unicyclic digraph in $\mathcal{U}(n)$ with maximum $\alpha$ trace norm.
We need the following singular value inequality for our main result \cite{ds}.\\
\begin{theorem}\cite{ds} Let $A$ and $B$ be square matrices of order $n$. Then $$\|A+B\|_*\le \|A\|_*+\|B\|_*.$$
Moreover, equality holds if and only if there exists an orthogonal matrix $P$ such that both $PA$ and $PB$ are positive semidefinite matrices.
\end{theorem}
\section{\bf Main Results}\label{sec2}
We first find variation in $\alpha$ trace norm of a digraph when an arc is deleted.
\begin{theorem}
Let $D$ be a digraph of order $n$ and with arc set $\mathcal{A(D)}$ and let $uv\in \mathcal{A(D)}$. Then\\$$\|{D_{\alpha}}\|_* \le \|{(D-uv)_{\alpha}}\|_{*}+\sqrt{2\alpha^2-2\alpha+1}$$\\
with equality if and only if $\alpha=0$ and $d^{+}(u)=d^{-}(v)=1$.
\end{theorem}
{\bf Proof.} We label the vertices such that $v$ is the first vertex, $u$ the second one. With this labelling, $A_{\alpha}(D)$ can be written as
\[ A_{\alpha}(D)=\begin{bmatrix} 
\alpha d_{v}^{+} ~& (1-\alpha) \delta ~& (1-\alpha) x^{T}\\
(1-\alpha) ~& \alpha d_u^{+}~&(1-\alpha) y^{T}\\
(1-\alpha) z ~& (1-\alpha)w~&E_{\alpha}
\end{bmatrix},
\] 
\vskip 3mm where $\delta \in\{0,1\}$, $x,y,z,w \in \mathbb{R}^{n-2}$ and $E_{\alpha}\in M_{n-2}$.\\
\vskip 3mm 
Also,\[ A_{\alpha}(D-uv)=\begin{bmatrix} 
\alpha d_{v}^{+} ~& (1-\alpha) \delta ~& (1-\alpha) x^{T}\\
0~& \alpha (d_u^{+}-1)~&(1-\alpha) y^{T}\\
(1-\alpha) z ~& (1-\alpha)w~&E_{\alpha}
\end{bmatrix}.
\] 
Therefore,\begin{align*}A_{\alpha}=A_{\alpha}(D-uv)+B_{\alpha},\end{align*}
where \[ B_{\alpha}=\begin{bmatrix} 
0 ~& 0~ & 0_{1\times (n-2)}\\
(1-\alpha) ~& \alpha ~&  0_{1\times (n-2)}\\
 0_{(n-2) \times 1}~ &  0_{(n-2) \times 1}~& 0_{(n-2) \times (n-2)}
\end{bmatrix}.
\]
\vskip 3mm 
Here we note that the singular values of $B_{\alpha}$ are $0^{[n-1]}$ and $\sqrt{2\alpha^{2}-2\alpha+1}$.
\vskip 3mm 
Now by Theorem $1.1$, \begin{align*}\|D_{\alpha}\|_{*}=\| {A_{\alpha}}(D)\|_{*}&=\| {A_{\alpha}}(D-uv)+B_{\alpha}\|_{*}\\
&\le \| {A_{\alpha}}(D-uv) \|_{*}+\|B_{\alpha}\|_{*}\\
&=  \| {(D-uv)}_{\alpha} \|_{*}+\sqrt{2\alpha^{2}-2\alpha+1},
\end{align*}
with equality if and only if there exists an orthogonal matrix $P$ such that $PA_{\alpha}(D-uv)$ and $PB_{\alpha}$ are both positive semidefinite matrices. Partition $P$ according to $A_{\alpha}(D-uv)$ as\\
\[ P=\begin{bmatrix} 
P_{11} ~& P_{12}~ & P_{13}\\
P_{21}~& P_{22}~&P_{23}\\
P_{31}~ & P_{32}~&P_{33}
\end{bmatrix}
\] 
\vskip 3mm 
Now,
\begin{align*}
P  B_{\alpha} &= 
\begin{bmatrix} 
P_{11} ~& P_{12}~ & P_{13} \\
P_{21}~ & P_{22}~ & P_{23} \\
P_{31}~ & P_{32}~ & P_{33} 
\end{bmatrix}
\begin{bmatrix} 
0~ & 0 ~& 0_{1\times (n-2)} \\
1-\alpha ~& \alpha ~& 0_{1\times (n-2)} \\
0_{(n-2) \times 1}~ & 0_{(n-2) \times 1}~ & 0_{(n-2) \times (n-2)} 
\end{bmatrix} \\
&=
\begin{bmatrix} 
(1-\alpha)P_{12}~ & \alpha P_{12} ~& 0_{1\times (n-2)}\\
(1-\alpha)P_{22}~ & \alpha P_{22} ~& 0_{1\times (n-2)}\\
(1-\alpha)P_{32} ~& \alpha P_{32} ~& 0_{(n-2) \times (n-2)} 
\end{bmatrix}
\end{align*}
\vskip 3mm 
Using the properties of positive semidefinite matrices (see Observation $7.1.10$ in \cite{hj}), we have $P_{32}=0$. Also, $P_{12}\ge 0$ and
$(1-\alpha)P_{22}=\alpha P_{12}$ and so $P_{22} \ge 0$.
\vskip 3mm 
Further using the properties of orthogonal matrix, we have $$P_{12}^2+P_{22}^2=1.$$
Two possibilities arise as $P_{12},P_{22} \ge 0$\\
$(i)$ $P_{12}=0$  and $P_{22}=1$.
$(ii$) $P_{12}=1, P_{22}=0$.\\
If $P_{12}=0$, then $(1-\alpha)P_{22}=0$ gives $P_{22}=0$. This means $P$ has a column of zeroes, which is not true.
So, $P_{12}=1$  and $P_{22}=0$.
Now, $P_{22}=0$ gives $\alpha P_{12}=0$, which gives $\alpha=0$ and $P_{12}=1$. Now, with  $ \alpha=0$, the matrix $PB_{\alpha}$ takes the following form\\
\[ PB_{\alpha}=\begin{bmatrix} 
1 ~& 0~ & 0_{1\times (n-2)}\\
0~ & 0 ~&0_{1\times (n-2)}\\
0_{(n-2) \times 1} ~& 0_{(n-2) \times 1}~&0_{(n-2) \times (n-2)} 
\end{bmatrix},
\] \\where 
\[ P=\begin{bmatrix} 
0 ~& 1~ & 0\\
P_{21}~ & 0~ & P_{23}\\
P_{31}~ & 0_{(n-2) \times 1}~&P_{33} 
\end{bmatrix}
\]\\by using orthonormal property of rows of $P$.\\
With $\alpha =0$\\ \[ A_{\alpha}(D-uv)=A_{0}(D-uv)=\begin{bmatrix} 
0 ~&  \delta~ &  x^{T}\\
0~& 0~& y^{T}\\
 z~ & w~&E_{0}
\end{bmatrix}.
\] \\
Now $P^{T}P=I_n$ gives\\
\[ \begin{bmatrix} 
0 ~& P_{21}^{T} ~& P_{31}^{T} \\
1 ~& 0 ~& 0_{1\times (n-2)} \\
0~& P_{23}^{T}~ & P_{33}^{T} 
\end{bmatrix} \begin{bmatrix} 
0 ~& 1 ~& 0\\
P_{21} ~& 0~ & P_{23}\\
P_{31} ~& 0_{(n-2) \times 1}~&P_{33} 
\end{bmatrix}
=I_n\]\\
or \[ \begin{bmatrix} 
P_{21}^{T}P_{21}+P_{31}^{T}P_{31} ~&  0 ~&  P_{21}^{T}P_{23}+P_{31}^{T}P_{33}\\
0~& 1~& 0\\
P_{23}^{T}P_{21}+P_{33}^{T}P_{31} ~ & 0 ~&P_{23}^{T}P_{23}+P_{33}^{T}P_{33}
\end{bmatrix}
= I_n.\]\\
This gives $P_{23}^{T}P_{23}+P_{33}^{T}P_{33}=I_{n-2}$.\\
Also, $PA_{\alpha}(D-uv)$ is a positive semidefinite, we have\\
\begin{align*}
PA_{\alpha}(D - uv) &= 
\begin{bmatrix} 
0 ~& 1 ~& 0 \\
P_{21}~ & 0~ & P_{23} \\
P_{31}~ & 0_{(n-2) \times 1} ~& P_{33} 
\end{bmatrix}
\begin{bmatrix} 
0 ~& \delta ~& x^{T} \\
0~ & 0 ~& y^{T} \\
z ~& w ~& E_{0} 
\end{bmatrix} \\
&= 
\begin{bmatrix} 
0 ~& 0 ~& y^{T} \\
P_{23}z ~& P_{21}\delta + P_{23}w ~& P_{21}x^{T} + P_{23}E_{0} \\
P_{33}z ~& P_{31}\delta + P_{33}w ~& P_{31}x^{T} + P_{33}E_{0} 
\end{bmatrix}.
\end{align*}

Now, using observation $7.1.10$ in \cite{hj}, we see $y^{T}=0$ which gives $d^{+}(u)=1$ \\
\vskip 3mm 
Also,
\begin{align*}
P_{23}z=0,P_{33}z&=0.\\
\mbox{This gives}~~~zz^{T}&=z^{T}(P_{23}^{T}P_{23}+P_{33}^{T}P_{33})z\\
 &=0\\
\implies z&=0~~~ \mbox{and ~so}~~~ d^{-}(v)=1.
\end{align*}
Conversely, if $d^{-}(v)=1=d^{+}(u)$ and $\alpha=0$. Then $y=0_{(n-2)\times 1}$ and $z=0_{(n-2)\times 1}$ respectively. Let $I_n(12)$ denote the permutation matrix obtained from $I_n$ by interchanging first two rows. Then the matrices $I_n(12)A_0(D)$ and $I_n(12)A_0(D-uv)$ are given by

\[ I_n(12)A_{0}(D)=\begin{bmatrix} 
1 ~& 0 ~& 0\\
0 ~& \delta~&x^{T}\\
0 ~& w~&E_{0}
\end{bmatrix}
\] 
and 
\[ I_n(12)A_{0}(D-uv)=\begin{bmatrix} 
0 ~& 0 ~& 0\\
0 ~& \delta~&x^{T}\\
0 ~& w~&E_{0}
\end{bmatrix}.
\] 
As $\|I_n(12)A_{0}(D)\|_*=\|A_{0}(D)\|_*$ and $\|I_n(12)A_{0}(D-uv)\|_*=\|A_{0}(D-uv)\|_*$, we see $\|D_0\|_*=\|(D-uv)_0\|_*+1$.\qed
A vertex $v$ in a digraph $D$ is said to be a leaf vertex if $d^{+}(v)+d^{-}(v)=1$. The next result as a consequence of Theorem $2.1$, gives the variation in $\alpha$ trace norm of a digraph by deletion of a leaf vertex.
 \begin{corollary}
 Let $D$ be a digraph with leaf vertex $u$ and $v$ the unique in neighbour or out neighbour of $u$. Then $\|D_{\alpha}\|_*\le \|(D-u)_{\alpha}\|_*+\sqrt{2{\alpha}^2-2\alpha+1}.$\\
 Moreover, equality holds if and only if $\alpha=0$ and either $d^+(u)=d^-(v)=1$ or $d^{-}(u)=d^{+}(v)=1$.
 \end{corollary}
{\bf Proof.} Removing the leaf vertex removes the arc between $u$ and $v$ in $D$. Recall that $\alpha$-trace norm remains unchanged by adding or removing isolated vertices, the result follows by Theorem $2.1$.\qed
A vertex $v$ in a digraph $D$ is said to be a nonleaf if $d^{+}(v)+d^{-}(v)\ge 2$. We next determine the variation in $\alpha$ trace norm of a digraph by deleting a nonleaf.
\begin{corollary}
Let $D$ be a digraph with $u$ a nonleaf vertex with indegree $d^{-}(u)$ and outdegree $d^{+}(u)$ so that total degree of $u$ is $d(u)=d^{-}(u)+d^{+}(u)$. Then $$\|D_{\alpha}\|_*\le \|(D-u)_{\alpha}\|_*+d(u)\sqrt{2{\alpha}^2-2\alpha+1},$$
with equality if and only if $\alpha=0$, $d^{+}(u)=1=d^{-}(u)$ and $d^{+}(w_1)=1$, $d^{-}(w_2)=1$, where $w_1$ is inneighbour of $u$ and $w_2$ is outneighbour of $u$.
\end{corollary}
{\bf Proof.} Removing nonleaf vertex $u$, removes all arcs between $u$ and its neighbours. The result now follows by Theorem $2.1$.\qed

Here we note that if we put $\alpha=0$ in Theorem $2.1$ and Corollary $2.2$, we respectively recover Theorem $3.5$ in \cite{mr2} and Theorem $2.7$ in \cite{apr}. In the next result, we determine the unique oriented tree in $\mathcal{T}(n)$ with maximum $\alpha$ trace norm.
\begin{corollary}
Let $T\in\mathcal{T}(n)$. Then 
\begin{eqnarray*}
\|T_{\alpha}\|_*\le (n-1)\sqrt{2{\alpha}^2-2\alpha+1}.
\end{eqnarray*}
 Moreover, equality holds if and only if $\alpha=0$ and $T=\overrightarrow{P_n}$.
\end{corollary}
{\bf Proof.} Let $T\in\mathcal{T}(n)$. Then $T$ has $n-1$ arcs and by Theorem $2.1$,  
\begin{eqnarray}
\|T_{\alpha}\|_*\le (n-1)\sqrt{2{\alpha}^2-2\alpha+1}.
\end{eqnarray}
It is easy to verify that in interval $[0,1)$, the function $\sqrt{2{\alpha}^2-2{\alpha}+1}$ decreases on $[0,\frac{1}{2}]$ and increases on $[\frac{1}{2},1)$ with maximum at $\alpha=0$ and minimum at $\alpha=\frac{1}{2}$. Assume that the equality holds in $(2.1)$. It is clear from the equality case in Theorem $2.1$ that $T$ cannot contain a vertex with indegree or outdegree greater or equal to 2. Hence only possibility is $\alpha=0$ and $T=\overrightarrow{P_n}$.\qed
Here we note that if we take $\alpha=0$ in Corollary $2.4$, we recover the Theorem $2.8$ in \cite{apr}. In the next result, we determine the unique digraph $U\in\mathcal{U}(n)$ with maximum $\alpha$-trace norm. We skip the proof as it is similar to Corollary $2.4$.\\
\begin{corollary}
Let $U\in\mathcal{U}(n)$. Then $\|U_{\alpha}\|_*\le n \sqrt{2{\alpha}^2-2\alpha+1}$. Equality holds if and only if $\alpha=0$ and $U=\overrightarrow{C_{n}}$.
\end{corollary}
{Acknowledgements} The research of Mushtaq A. Bhat is supported by SERB-DST grant with File No. MTR/2023/000201 and by NBHM project number NBHM/02011/20/2022. The research of Peer Abdul Manan is supported by CSIR, New Delhi, India with CSIR-HRDG Ref. No: Jan-Feb/06/21(i)EU-V. 

\end{document}